\documentclass[preprint,times]{elsarticle}

\usepackage{graphicx}%
\usepackage{multirow}%
\usepackage{amsmath,amssymb,amsfonts}%
\usepackage{amsthm}%
\usepackage{mathrsfs}%
\usepackage[title]{appendix}%
\usepackage{xcolor}%
\usepackage{textcomp}%
\usepackage{manyfoot}%
\usepackage{booktabs}%
\usepackage{algorithm}%
\usepackage{algorithmicx}%
\usepackage{algpseudocode}
\usepackage{color}%
\usepackage{caption}
\usepackage{subcaption}
\usepackage{listings}%
\newtheorem{theorem}{Theorem}
\newtheorem{lemma}{Lemma}%
\newtheorem{proposition}{Proposition}%
\newtheorem{remark}{Remark}%

\begin{document}
\begin{frontmatter}
\title{Time integration of dissipative stochastic PDEs}

\author[1]{Helena Biščević}
\author[2]{Raffaele D'Ambrosio}

\address[1]{Gran Sasso Science Institute, L'Aquila, Italy}
\address[2]{Department of Information Engineering and Computer Science and Mathematics, University of L'Aquila, Italy}

\begin{abstract}
The paper is focused on the numerical solution of stochastic reaction-diffusion problems. A special attention is addressed to the conservation of mean-square dissipativity in the time integration of the spatially discretized problem, obtained by means of finite differences. The analysis highlights the conservative ability of stochastic $\theta$-methods and stochastic $\theta$-IMEX methods, emphasizing the roles of spatial and temporal stepsizes. A selection of numerical experiments is provided, confirming the theoretical expectations.
\end{abstract}
\begin{keyword}
Stochastic partial differential equations, stochastic $\theta$-methods, stochastic $\theta$-IMEX methods, mean-square dissipativity
\MSC[2020] 65C30 \sep 60H10 \sep 60H15
\end{keyword}
\end{frontmatter}

\section{Introduction}
We focus our attention on the numerical solution of stochastic partial differential equations (SPDEs), whose meaningful role in modeling a wide range of real-life phenomena has clearly been highlighted by the existing literature, which also analyzes their properties of well posedness, stability, regularity of the solutions (see, for instance, \cite{bax,daprato,holden,lui,loto,weinan,spde_kuksin} and references therein). As a consequence, the numerical integration of SPDEs has been the core topic of many recent key contributions, mostly based on the establishment of numerical schemes and the analysis of their accuracy (we can refer, for instance, to the non-exhaustive selection \cite{karnia,kruse,mt04,hongInv} and references therein).

A specific class of problems treated in this paper is given by stochastically forced dissipative nonlinear partial differential equations, formalized in terms of the stochastic reaction-diffusion problem
\begin{equation}\label{generalSPDE}
    {\rm d}u(x,t)=\left(-\mathcal{A}u(x,t)+\mathcal{R}(u(x,t))\right){\rm d}t +g(u(x,t)) {\rm d} W(x,t),
\end{equation}
where $t\geq 0$, equipped by the initial value $u(x,0)=u_0$. The problem is assumed to be 1D with specified boundary conditions. The solution to this problem belongs to a separable Hilbert space $\mathbb{H}$ endowed with the inner product $\langle \cdot,\cdot \rangle_{\mathbb{H}}$ and its existence and uniqueness is proved in \cite{daprato2}. The operator $\mathcal{A}$ is linear and self-adjoint on a domain $D(A)\subset\mathbb{H}$ and a complete orthonormal system of eigenvectors $\{e_k(x), \ k\in\mathbb{N}\}$ is forming a basis in $L^2(D(A))$. The nonlinear functions $\mathcal{R}$ and $g$ in \eqref{generalSPDE} are respectively defined in the domains $D(R),D(g)\subset\mathbb{H}$ and assume values in $\mathbb{H}$. 

The noisy term in \eqref{generalSPDE} is generally multiplicative even though, in the remainder of the paper, the analysis is carried out separately for additive and multiplicative noises, for methodological purposes. It is governed by the cylindrical Wiener process
$$W(x,t)=\sum_{k\in\mathbb{N}} \sigma_k(t)\omega_k(t)e_k(x),$$
where $\sigma_k(t)\in\mathbb{R}$, $t\geq 0$, $\{\omega_k(t), \ t\geq 0\}_{k\in\mathbb{N}}$ is a set of independent standard Wiener processes. 

It is worth observing that Equation \eqref{generalSPDE} is a stochastic perturbation of a deterministic reaction-diffusion problem, via multiplicative noise. The dissipative character of the problem is studied, for instance, in \cite{daprato,mhairer,spde_kuksin,weinan}, mostly by means of invariant measure arguments. For such problems, the present paper aims to address a key question about the ability of numerical methods to preserve the dissipative character, with respect to selected time integrators applied to the spatial discretization, obtained by means of finite differences. Let us point out that, for genuine SDE-based problems, the numerical conservation of the dissipative character is generally translated into proper stepsize restrictions that a prescribed method has to satisfy \cite{bdd2025,bd2021,dambook,DAMBROSIO2021105671,dadigioOptim}. Here, we aim to specify that, if a system of SDEs arises from the spatial discretization of a dissipative SPDE, the contractivity is either preserved or not, depending on the method and the value of $\theta$. 

The paper is organized as follows: Section 2 provides the spatial discretization of the problem under consideration; Section 3 focuses on the analysis of dissipativity in mean-square sense for $\theta$-methods, both for additive and multiplicative noises; an analogous investigation is object of Section 4; Section 5 provides the monotonicity analysis under one-sided Lipschitz assumption for the drift of the spatially discretized problem; the numerical evidence arising from selected experiments is then presented in Section 6; concluding remarks and future issues regarding the prosecution of this research are highlighted in Section 7.

\section{Spatial discretization} 
Let us compute the semi-discrete version of Equation \eqref{generalSPDE}, where the spatial variable is fully discretized by means of finite differences. In order to proceed towards this direction, we provide a uniform partition of the interval $[x_{\ell},x_{r}]$, given by 
$${\cal I}_{\Delta x}=\left\{x_i=x_{\ell}+i\Delta x, \quad i=0,1,\ldots,d+1, \quad \Delta x=\frac{x_r-x_{\ell}}{d+2}\right\}.$$
In order to discretize the differential operator ${\cal A}$ by means of centered finite differences, we define the matrix 
$$A=\frac1{\Delta x^p}B,$$
where $B$ is a band matrix of length $p+1$, $p$ being the order of the employed finite difference \cite{bookfindiff}. Moreover, we introduce the vectors 
$$U(t)=\left[u(x_i,t)\right]_{i=0}^{d+1}, \quad W(t)=\left[W(x_i,t)\right]_{i=0}^{d+1}.$$ 
Correspondingly, the general form of the spatial semi-discretization of Equation \eqref{generalSPDE} is given by the following system of SDEs
\begin{equation}\label{spatially_discr}
    {\rm d}U(t) \approx \left(-A U(t) + R(U(t))\right){\rm d}t + G(U(t)) {\rm d} W(t).
\end{equation}
where
$$R(U(t))=\begin{bmatrix}
\mathcal{R}(u(x_0,t))\\[2mm]
\mathcal{R}(u(x_1,t))\\[2mm]
\vdots\\[2mm]
\mathcal{R}(u(x_{d+1},t))
\end{bmatrix}$$
and
$$G(U(t))=\begin{bmatrix}
g(u(x_0,t))\\[2mm]
g(u(x_1,t))\\[2mm]
\vdots\\[2mm]
g(u(x_{d+1},t))
\end{bmatrix}.$$

In the remainder, we proceed according to the following assumptions.

\medskip

\noindent{\bf Assumption 1}.\\ The function $R(U(t))$ is of class $C^1(\mathbb{R}^{d+2})$.

\medskip

\noindent{\bf Assumption 2}.\\ The function $G(U(t))$ is locally Lipschitz continuous, i.e., for any $U^*\in\mathbb{R}^{d+2}$ there exist $L_g>0$ and $\rho>0$ such that, 
$$\|G(U)-G(U^*)\|\leq L_g\|U-U^*\|,$$
for any $U\in\mathbb{R}^{d+2}$ such that $\|U-U^*\|\leq\rho$.

\medskip 


The purpose of the paper is the numerical conservation of dissipativity, supposing that the original problem is dissipative in mean-square sense. Such a property will then be checked on the case-by-case basis in the numerical experiments.  

\section{Mean-square dissipativity of $\theta$-methods}
In the following section, we will present results regarding additive and multiplicative noises, using $\theta$-Maruyama method for time discretization of problem \eqref{spatially_discr}. This choice is emerging from its easy implementation for this type of equations as well as its wide usage in applications. We would like to see how much we can exploit simple methods in this framework, as by looking at the previous results done on SDEs, it showed to be reasonable to use \cite{DAMBROSIO2021105671,bdd2025}.

Referring to the discretized domain 
\begin{equation}\label{ddom}
{\cal I}_{\Delta t}=\left\{t_n=n\Delta t, \quad n=0,1,\ldots,N, \quad \Delta t=\frac{T}{N}\right\},
\end{equation}
the aforementioned method, applied to a general system of SDEs having the form
$${\rm d}X(t)=f(X(t)){\rm d}t+g(X(t)){\rm d}W(t),$$
reads
\begin{equation}\label{teta}
    X_{n+1} = X_n + (1-\theta)\Delta t f (X_n) +  \theta \Delta t f(X_{n+1}) + g(X_n)\Delta W_n,
\end{equation}
where $\theta \in [0,1]$, $X_n$ is the approximate value for $X(t_n)$ and the Wiener increment $\Delta W_n=W(t_{n+1})-W(t_n)$ is distributed as a Gaussian random variable with zero mean and variance $\Delta t$.

In the proofs of contractivity results from this section on, we will be using the estimates on the matrices contained in the following lemma.

\medskip

\begin{lemma}\label{estimates}
    Let $D(\theta, \Delta t)$ and $C(\theta, \Delta t)$ be defined as 
    $$\begin{aligned}
        D(\theta, \Delta t)&=\left(I + \theta \Delta t A\right)^{-1},\\[1mm]
        C(\theta, \Delta t)&=D(\theta, \Delta t)\left( I - (1-\theta) \Delta t A\right),
    \end{aligned}$$
    where $I\in\mathbb{R}^{d+2}$ is the identity matrix, $A$ is a Toeplitz matrix, $\theta \in [0,1]$ and $\Delta t >0$. Then, the following estimates hold:

$$\begin{aligned}
     \text{(i)} & \quad \left\|D(\theta, \Delta t)\right\|\leq\frac{1}{1-\theta\Delta t\|A\|}, \\[2mm]
     \text{(ii)} & \quad \left\|C(\theta, \Delta t)\right\|\leq \frac{1+(1-\theta)\Delta t\|A\|}{1-\theta\Delta t\|A\|}.
\end{aligned}$$

\end{lemma}
\begin{proof}
Estimate $(i)$ derives from reverse triangle inequality. Regarding $(ii)$, the proof follows using $(i)$ and the estimate $\|I - (1-\theta) \Delta t A\|\leq 1+ (1-\theta) \Delta t\|A\|$.
\end{proof}

\subsection{Analysis in case of additive noise}\label{thetaAdd}

We first address the analysis of \eqref{teta}, in presence of additive noise. The following result highlights properties of the mean-square deviation between two solutions that, to some extent, imitate the underlying deterministic case, as it will be clarified in the proof. 

\begin{theorem}
    For the spatially discretized problem \eqref{spatially_discr}, with two distinct initial values $U_0$ and $Y_0$, under Assumption 1,  the numerical dynamics generated by $\theta$-Maruyama method \eqref{teta} satisfies the inequality
    \begin{equation}
    \label{contrCondthadd}
        \mathbb{E}\left\| U_{n} - Y_{n} \right\|^2 \le \alpha(\theta,\Delta t)\mathbb{E}\left\| U_{0} - Y_{0} \right\|^2,
    \end{equation}
    where 
    \begin{equation}\label{coeff1}
    \alpha(\theta,\Delta t)=\left( \frac{1+(1-\theta) \Delta t (\left\|A\right\| + M)}{1-\theta \Delta t (\left\|A\right\| + M)} \right)^2,
    \end{equation}
    with
    \begin{equation}\label{M}
    M=\sup_{t \in [0,T]}{\mathbb{E}\left[\left\|R'(U(t))\right\|^2\right]}.
    \end{equation}
\end{theorem}
\begin{proof}
    Applying $\theta$-Maruyama method \eqref{teta} for time discretization of \eqref{spatially_discr}, with initial value $U_0$, we obtain
    \begin{equation}
        U^{n+1} = U^{n} + (1-\theta) \Delta t [-AU^n + R(U^n)] + \theta \Delta t [-AU^{n+1} + R(U^{n+1})] + {\epsilon} \Delta W_{n}
    \end{equation}
    In order to study mean-square contractivity, we proceed similarly for problem \eqref{spatially_discr} with initial value $Y_0$, i.e.,
    \begin{equation}
        Y^{n+1} = Y^{n} + (1-\theta) \Delta t [-AY^n + R(Y^n)] + \theta \Delta t [-AY^{n+1} + R(Y^{n+1})] + {\epsilon}\Delta W_{n}
    \end{equation}
    Then, by denoting $Z(t)=U(t)-Y(t)$ and proceeding with side-by-side subtraction, we have 
    %
    \begin{equation}
    \begin{aligned}
    Z^{n+1} = C(\theta,\Delta t) Z^n &+ (1-\theta) \Delta t D(\theta,\Delta t)\left(R(U^n)-R(Y^n)\right)\\[1mm]
    &+ \theta \Delta t D(\theta,\Delta t)\left(R(U^{n+1}) - R(Y^{n+1})\right),
    \end{aligned}
    \end{equation}
    where $D(\theta, \Delta t)$ and $C(\theta, \Delta t)$ are defined as in Lemma \ref{estimates}.
    Passing to the norm and by definition of $M$, given in \eqref{M}, we get 
    \begin{equation}
        \left\|Z^{n+1}\right\| \leq \left\|C(\theta, \Delta t)\right\| \left\|Z^{n}\right\| + (1-\theta) \Delta t \left\| D(\theta, \Delta t) \right\|M \left\|Z^{n}\right\| + \theta \Delta t \left\|D(\theta, \Delta t)\right\| M \left\|Z^{n+1}\right\|.
    \end{equation}
    Using estimates from Lemma \ref{estimates}, we obtain
    \begin{equation}
         \left\|Z^{n+1}\right\| \leq \frac{1+(1-\theta) \Delta t (\left\|A\right\| + M)}{1-\theta \Delta t (\left\|A\right\| + M)} \left\|Z^{n}\right\|.
    \end{equation}
    Finally, squaring and passing to the expectation leads to the thesis.
\end{proof}

Supposing that the mean-square deviation between two solutions of the problem is monotonically decreasing, let us study the contractivity condition $\alpha(\theta,\Delta t)<1$ for the method, i.e., 
\begin{equation}\label{eqpar}
    (1 - 2\theta)\Delta t(\|A\|+M)<-2.
\end{equation}
We distinguish the following two cases:
\begin{itemize}
\item for $(1 - 2\theta)>0$, i.e., when $\theta<\frac12$, Equation \eqref{eqpar} becomes 
$$\Delta t<-\frac2{(1 - 2\theta)(\|A\|+M)}.$$
Since the quantity on the right-hand side is always negative, in this case $\theta$-methods do not preserve the contractivity character of the spatially discretized problem;
\item for $(1 - 2\theta)<0$, i.e., when $\theta>\frac12$, since $A=\frac1{\Delta x^p}B$, we can rewrite the inequality given by Equation \eqref{eqpar} in terms of a relation between the two stepsizes, i.e., 
$$\frac{\Delta t}{\Delta x^p}>\frac1{\|B\|}\left(- M\Delta t-\frac2{(1-2\theta)}\right).$$
We observe that the contractive behaviour of the numerical dynamics depends on the sign of its right-hand side. In particular, when it is negative, the method provides a contractive dynamics, for any value of $\Delta t$ and $\Delta x$, without restrictions. Instead, when it is positive, the method does not show any contractive behaviour.
\end{itemize}

\subsection{Analysis in case of multiplicative noise}\label{thetaMul}
After analysing the additive noise, we proceed with considering a multiplicative noise. The proof of this result, compared to the proof with additive noise, will require additional hypothesis on the diffusion function as well as further handling of the noise term. However, in the following, we will show that the result is achievable in an intuitive and straightforward way. 

\begin{theorem}
    For the spatially discretized problem \eqref{spatially_discr} with two distinct initial values $U_0$ and $Y_0$, under Assumptions 1 and 2, the numerical dynamics generated by $\theta$-Maruyama method \eqref{teta} satisfies the inequality
    \begin{equation}
    \label{contrCondthmul}
        \mathbb{E}\left\| U_{n} - Y_{n} \right\|^2 \le \alpha(\theta,\Delta t)\mathbb{E}\left\| U_{0} - Y_{0} \right\|^2,
    \end{equation}
    where 
    \begin{equation}\label{coeff2}
    \alpha(\theta,\Delta t)=\frac{(1+(1-\theta) \Delta t (\left\|A\right\| + M))^2 + L_g^2 \Delta t}{(1-\theta \Delta t (\left\|A\right\| + M))^2},
    \end{equation}
    with 
\begin{equation}\label{M2}
    M=\sup_{t \in [0,T]}{\mathbb{E}\left\|R'(U(t))\right\|}.
    \end{equation}
\end{theorem}
\begin{proof}
    Applying $\theta$-Maruyama method \eqref{teta} for time discretization of \eqref{spatially_discr} to compute two distinct solutions, with initial values $U_0$ and $Y_0$, we have
    $$
        U^{n+1} = U^{n} + (1-\theta) \Delta t [-AU^n + R(U^n)] + \theta \Delta t [-AU^{n+1} + R(U^{n+1})] + g(U^n) \Delta W_{n},
    $$
    $$
        Y^{n+1} = Y^{n} + (1-\theta) \Delta t [-AY^n + R(Y^n)] + \theta \Delta t [-AY^{n+1} + R(Y^{n+1})] + g(Y^n) \Delta W_{n}.
   $$
    Then, by denoting $Z(t)=U(t)-Y(t)$ and proceeding with side-by-side subtraction, we have  
    \begin{equation*}
    \begin{aligned}
    Z^{n+1} = C(\theta,\Delta t) Z^n &+ (1-\theta) \Delta t D(\theta,\Delta t)\left(R(U^n)-R(Y^n)\right)\\[1mm]
    &+ \theta \Delta t D(\theta,\Delta t)\left(R(U^{n+1}) - R(Y^{n+1})\right) + (g(U^n)-g(Y^n)) \Delta W_{n},
    \end{aligned}
    \end{equation*}
    where the matrices $D(\theta, \Delta t)$ and $C(\theta, \Delta t)$ are defined by Lemma \ref{estimates}, as before.
    Using the estimates from Lemma \ref{estimates}, Assumption 2 and the definition of $M$, given by \eqref{M2}, we obtain
    \begin{equation}
         \frac{1-\theta \Delta t (M + \left\|A\right\|)}{1-\theta \Delta t \left\|A\right\|}\left\|Z^{n+1}\right\| \leq \beta(\theta, \Delta t) \left\|Z^{n}\right\|
    \end{equation}
    where
    \begin{equation}
        \beta(\theta, \Delta t) = \frac{1+(1-\theta) \Delta t (\left\|A\right\| + M) + L_g \left\| \Delta W_n \right\|}{1-\theta \Delta t \left\|A\right\|}
    \end{equation}
    By squaring and passing to the expectation, we have
    \begin{equation}
    \begin{aligned}
        \mathbb{E} \beta(\theta, \Delta t)^2 &= \frac{1 + (1-\theta)^2 \Delta t ^2 (\left\|A\right\| + M)^2 + L_g^2\Delta t + 2(1-\theta) \Delta t (\left\|A\right\| + M)}{(1-\theta \Delta t \left\|A\right\|)^2}.
    \end{aligned}
    \end{equation}
    After discarding the terms that are zero in expectation, we get the result.
\end{proof}

Supposing that the mean-square deviation between two solutions of the problem is monotonically decreasing, let us study the contractivity condition $\alpha(\theta,\Delta t)<1$ for the method, i.e.,
\begin{equation}\label{eqpar2multheta}
    \Delta t (2(\|A\| + M)+\Delta t (1 - 2\theta)(\|A\| + M)^2+L_g^2)<0.
\end{equation}


We observe that
\begin{itemize}
\item for $(1 - 2\theta)>0$, i.e., when $\theta<\frac12$, we have
\begin{equation*}
    0<\Delta t < -\frac{2(\|A\| + M)+L_g^2}{(1 - 2\theta)(\|A\| + M)^2}.
\end{equation*}
The inequality \eqref{eqpar2multheta} is not fulfilled for any value of $\Delta t$ and, as a direct consequence, $\theta$-methods do not preserve the contractivity in this case;
\item for $(1 - 2\theta)<0$, i.e., when $\theta>\frac12$, inequality \eqref{eqpar2multheta} reduces to 
\begin{equation*}
    \Delta t >-\frac{2(\frac1{\Delta x^p}\|B\| + M)+L_g^2}{(1-2\theta)(\frac1{\Delta x^p}\|B\| + M)^2}.
\end{equation*}
We observe that, for $\Delta x$ tending to 0, the right-hand side of last inequality tends to 0, so it is always fulfilled for any choice of $\Delta t$. In other terms, $\theta$-methods for $\theta>\frac12$ are unconditionally contractive in time, for reasonably small values of $\Delta x$. 
\end{itemize}


\begin{remark}
We draw attention to the case of $\theta=\frac12$, that is not covered in the analysis above, both for additive and multiplicative noises. This case would need a separate investigation, relying on alternative arguments. However, the fact that the existing analysis covers all the other values of $\theta>\frac12$ is leading to a satisfactory behaviour, in terms of structure preservation, makes this further investigation dispensable.  

We also observe that the case $\theta=\frac12$ is instead covered by the following analysis regarding $\theta$-IMEX methods.
\end{remark}

\section{Mean-square dissipativity of $\theta$-IMEX methods}

Methods of IMEX type allow us to tackle stiff problems, by dividing the part of the problem that needs to be treated implicitly, from the part that can be safely solved by explicit methods. Here we arrive to the omnipresent problem of choosing explicit and implicit methods in numerical analysis as both bring benefits, but also detriments in their further analysis. In our case, it is worth trying to apply the following approach: reaction is treated explicitly and diffusion of the deterministic part implicitly.

Discretizing the domain as in Equation \eqref{ddom}, $\theta$-IMEX methods applied to the general system of SDEs 
$${\rm d}X(t)=\bigl(\Lambda X(t)+f(X(t))\bigr){\rm d}t+g(X(t)){\rm d}W(t),$$
read
\begin{equation}\label{tetaIMEX}
    X_{n+1} = X_n + (1-\theta)\Delta t\Lambda X_n+\theta\Delta t\Lambda X_{n+1}+\theta \Delta t f(X_{n}) + g(X_n)\Delta W_n,
\end{equation}
where the matrix $\Lambda\in\mathbb{R}^{(d+2)\times (d+2)}$ arises from the discretization of the differential operator ${\cal A}$ in \eqref{generalSPDE}.

\subsection{Analysis in case of additive noise}

Proceeding in analogous way as in Section \ref{thetaAdd}, we first analyse the behaviour of $\theta$-IMEX methods \eqref{tetaIMEX} in presence of additive noise. \\

\begin{theorem}
    For the spatially discretized problem \eqref{spatially_discr}, with two distinct initial values $U_0$ and $Y_0$, under Assumption 1, the numerical dynamics generated by $\theta$-IMEX method satisfies the inequality
    \begin{equation}
    \label{contrCondIMadd}
        \mathbb{E}\left\| U_{n} - Y_{n} \right\|^2 \le \alpha(\theta,\Delta t)\mathbb{E}\left\| U_{0} - Y_{0} \right\|^2,
    \end{equation}
    where 
    \begin{equation}\label{coeff3}
    \alpha(\theta,\Delta t)=\left( \frac{1+(1-\theta) \Delta t \left\|A\right\| + \Delta t M}{1-\theta \Delta t \left\|A\right\|} \right)^2,
    \end{equation}
    with
    \begin{equation}\label{M3}
    M=\sup_{t \in [0,T]}{\mathbb{E}\left\|R'(U(t))\right\|}.
    \end{equation}
\end{theorem}
\begin{proof}
    Applying $\theta$-IMEX method \eqref{teta} for time discretization of \eqref{spatially_discr}, with initial value $U_0$, we obtain
    \begin{equation}
        U^{n+1} = U^{n} + (1-\theta) \Delta t (-AU^n) + \theta \Delta t (-AU^{n+1}) + \Delta t R(U^{n}) + {\epsilon} \Delta W_{n}
    \end{equation}
    In order to study mean-square contractivity, we proceed similarly for problem \eqref{spatially_discr} with initial value $Y_0$, i.e.,
    \begin{equation}
        Y^{n+1} = Y^{n} + (1-\theta) \Delta t (-AY^n) + \theta \Delta t (-AY^{n+1}) + \Delta t R(Y^n) + {\epsilon}\Delta W_{n}
    \end{equation}
    Then, by denoting $Z(t)=U(t)-Y(t)$ and proceeding with side-by-side subtraction, we have 
    \begin{equation}
    Z^{n+1} = C(\theta,\Delta t) Z^n + \Delta t D(\theta,\Delta t)\left(R(U^n)-R(Y^n)\right),
    \end{equation}
    where $D(\theta, \Delta t)$ and $C(\theta, \Delta t)$ are defined by Lemma \ref{estimates}. Passing to the norm, using Lemma \ref{estimates} and the definition of $M$, as in \eqref{M3}, 
    \begin{equation}
        \left\|Z^{n+1}\right\| \leq \left\|C(\theta, \Delta t)\right\| \left\|Z^{n}\right\| + \Delta t \left\| D(\theta, \Delta t) \right\|M \left\|Z^{n}\right\|.
    \end{equation}
    Using estimates from Proposition \ref{estimates}, we obtain
    \begin{equation}
         \left\|Z^{n+1}\right\| \leq \frac{1+(1-\theta) \Delta t \left\|A\right\| + \Delta t M}{1-\theta \Delta t \left\|A\right\|} \left\|Z^{n}\right\|.
    \end{equation}
    Finally, squaring and passing to the expectation leads to the thesis.
\end{proof}

Supposing that the mean-square deviation between two solutions of the problem is monotonically decreasing, let us study the contractivity condition $\alpha(\theta,\Delta t)<1$ for the method, i.e., 
\begin{equation}\label{eqpar3}
    \Delta t (2(1-\theta \Delta t\|A\|)+\Delta t (\|A\| + M))<0.
\end{equation}
Therefore, contractivity sums up to the following:
\begin{equation}\label{eqpar2bis}
    2 + \Delta t (M + (1-2\theta)\|A\|)<0.
\end{equation}
Hence, we can distinguish two main cases: 
\begin{itemize}
    \item for $(1 - 2\theta)>0$, i.e., when $\theta<\frac12$, contractivity condition translates into
    \begin{equation*}
        0<\Delta t < -\frac{2}{(1 - 2\theta)\|A\| + M}.
    \end{equation*}
    Analysing the quantity on the right-hand side, by taking $A=\frac1{\Delta x^p}B$, reveals that
    \begin{equation*}
        \lim_{\Delta x \to 0} 
        -\frac{2 \Delta x^p}{(1 - 2\theta)\|B\| + \Delta x^p M} = 0,
    \end{equation*}
    so the method does not preserve contractivity in this case;
    \item for $(1 - 2\theta)\leq 0$, i.e., when $\theta\geq\frac12$, the condition
the condition $M+(1-2\theta)\|A\|<0$ brings us to the spatial stepsize restriction
\begin{equation*}
        \Delta x^p < -\frac{(1-2\theta)\|B\|}{M}.
    \end{equation*}
    Let us observe that such a restriction can be rather severe for large enough values of $M$.
By replacing $\Delta x^p\|A\|=\|B\|$ in \eqref{eqpar2bis}, we have
$$2\Delta x^p+\Delta t \left(\Delta x^p + (1-2\theta)\|B\|\right)<0$$
and, in the limit for $\Delta x \to 0$, it reads as
$$\Delta t (1-2\theta)\|B\|<0.$$
Last inequality is clearly fulfilled for any $\Delta t$. In other terms, $\theta$-IMEX methods are unconditionally contractive only for infinitesimal values of $\Delta x$, making them unfeasible for time integration of dissipative SPDEs. 
\end{itemize}

\subsection{Analysis in case of multiplicative noise}

Similarly as in Section \ref{thetaMul}, we now analyse the behaviour of $\theta$-IMEX methods \eqref{tetaIMEX} in presence of multiplicative noise. \\

\begin{theorem}
    For the spatially discretized problem \eqref{spatially_discr}, with two distinct initial values $U_0$ and $Y_0$, under Assumptions 1 and 2,  the numerical dynamics generated by $\theta$-IMEX methods satisfies the inequality
    \begin{equation}
    \label{contrCondIMmul}
        \mathbb{E}\left\| U_{n} - Y_{n} \right\|^2 \le \alpha(\theta,\Delta t)\mathbb{E}\left\| U_{0} - Y_{0} \right\|^2,
    \end{equation}
    where 
    \begin{equation}\label{coeff4}
    \alpha(\theta,\Delta t)= \frac{(1+(1-\theta) \Delta t \left\|A\right\| + \Delta t M)^2 + L_g^2 \Delta t}{(1-\theta \Delta t \left\|A\right\|)^2},
    \end{equation}
    with
    \begin{equation}\label{M4}
    M=\sup_{t \in [0,T]}{\mathbb{E}\left\|R'(U(t))\right\|}.
    \end{equation}
\end{theorem}
\begin{proof}
    Applying $\theta$-IMEX method \eqref{teta} for time discretization of \eqref{spatially_discr}, with initial value $U_0$, we obtain
    \begin{equation}
        U^{n+1} = U^{n} + (1-\theta) \Delta t (-AU^n) + \theta \Delta t (-AU^{n+1}) + \Delta t R(U^{n}) + g(U^n) \Delta W_{n}
    \end{equation}
    In order to study mean-square contractivity, we proceed similarly for problem \eqref{spatially_discr} with initial value $Y_0$, i.e.,
    \begin{equation}
        Y^{n+1} = Y^{n} + (1-\theta) \Delta t (-AY^n) + \theta \Delta t (-AY^{n+1}) + \Delta t R(Y^n) + g(Y^n) \Delta W_{n}
    \end{equation}
    Then, by denoting $Z(t)=U(t)-Y(t)$ and proceeding with side-by-side subtraction, we have 
    \begin{equation}
    Z^{n+1} = C(\theta,\Delta t) Z^n + \Delta t D(\theta,\Delta t)\left(R(U^n)-R(Y^n)\right) + D(\theta,\Delta t) (g(U^n) - g(Y^n)) \Delta W_n
    \end{equation}
    where the matrices $D(\theta, \Delta t)$ and $C(\theta, \Delta t)$ are those defined by Lemma \ref{estimates}. Passing to the norm, using the definition of $M$ as in \eqref{M4} and Assumption 2 yields 
    \begin{equation}
        \left\|Z^{n+1}\right\| \leq \left\|C(\theta, \Delta t)\right\| \left\|Z^{n}\right\| + \Delta t \left\| D(\theta, \Delta t) \right\|M \left\|Z^{n}\right\| + \left\| D(\theta, \Delta t) \right\| L_g \left\|Z^{n}\right\| \left\| \Delta W_n \right\|.
    \end{equation}
    Using estimates from Lemma \ref{estimates}, we obtain

    \begin{equation}
         \left\|Z^{n+1}\right\| \leq \frac{1+(1-\theta) \Delta t \left\|A\right\| + \Delta t M + L_g \left\| \Delta W_n \right\|}{1-\theta \Delta t \left\|A\right\|} \left\|Z^{n}\right\|.
    \end{equation}
    Finally, squaring and passing to the expectation leads to the thesis.
\end{proof}

\medskip

Supposing that the mean-square deviation between two solutions of the problem is monotonically decreasing, let us study the contractivity condition $\alpha(\theta,\Delta t)<1$ for the method, i.e.,
\begin{equation}\label{eqpar2}
    \Delta t (\Delta t (M^2 + (1-2\theta)\|A\|^2 +2(1-\theta)M\|A\|) + 2M+2\|A\| + L_g^2)<0.
\end{equation}

Hence, we can distinguish two main cases: 
\begin{itemize}
    \item for $(1 - 2\theta)>0$, i.e., when $\theta<\frac12$, contractivity condition translates into
    \begin{equation*}
        0<\Delta t < -\frac{2(M+\|A\|)+L_g^2}{M^2 + (1 - 2\theta)\|A\|^2 + 2(1-\theta)M\|A\|}.
    \end{equation*}
    Analysing the quantity on the right-hand side, by taking $A=\frac1{\Delta x^p}B$, reveals that
    \begin{equation*}
        \lim_{\Delta x \to 0} 
        -\frac{2(M \Delta x^p + \|B\|) + L_g^2\Delta x^p}{\Delta x^{2p} M^2 + (1 - 2\theta)\|B\|^2 + 2(1-\theta)M\|B\|\Delta x^p}\Delta x^p = 0,
    \end{equation*}
    so the method does not preserve contractivity in this case;
    \item for $(1 - 2\theta)\leq 0$, i.e., when $\theta\geq\frac12$, the analysis of
    contractivity proceeds in analogous way as in the previous case. In particular, the condition $M^2 + (1 - 2\theta)\|A\|^2 + 2(1-\theta)M\|A\|<0$ yields the following restriction on $\Delta x$, by taking $A=\frac1{\Delta x^p}B$, i.e.,
    \begin{equation*}
        \Delta x^p \in \left(-\frac{\|B\|(1+2\theta)}{M}, -\frac{\|B\|}{M}\right),
    \end{equation*}
    which is impossible to satisfy. Therefore, we expect $\theta$-IMEX method to not be as efficient when attempting to preserve contractivity for stochastic reaction-diffusion PDE with multiplicative noise.
\end{itemize}








\section{Monotonicity analysis under one-sided Lipschitz drift}
The analysis so far has been based on the assumptions that $R$ is $C^1$ and $G$ is locally Lipschitz. We now aim to understand what happens when the drift and the diffusion of SDEs arising from the spatial discretization of \eqref{generalSPDE} satisfy the classical hypotheses of one-sided and global Lipschitz continuity, on drift and diffusion respectively, required in the existing literature (see, for instance, \cite{bdd2025a,h00,DAMBROSIO2021105671,dambook,r7}) to achieve the monotonicity of the mean-square deviation between two solutions. 

Accordingly, in this section, we proceed relying on the following assumptions related to Equation \ref{spatially_discr}.

\medskip

{\bf Assumption 3}.
The function $R(U(t))$ is one-sided Lipschitz continuous with constant $\mu\in\mathbb{R}$, i.e., $$\langle R(U) - R(V), U-V \rangle \leq \mu \|U-V\|^2.$$

{\bf Assumption 4}.
The function $G(U(t))$ is globally Lipschitz continuous with constant $L>0$, i.e.,
$$\|G(U)-G(V)\|^2\leq L\|U-V\|^2.$$

\medskip

Let us provide the following useful result. 

\medskip

\begin{proposition}
Let $R$ be one-sided Lipschitz with constant $\mu$. Then, the whole drift of \eqref{spatially_discr} is also one-sided Lipschitz with constant $\mu^{*}=\|A\|+\mu$.
\end{proposition}
\begin{proof}
The result follows from the direct verification of one-sided Lipschitz continuity for the whole drift of \eqref{spatially_discr}, i.e., for any $x,y\in\mathbb{R}^{d+2}$,
    \begin{align*}
        \langle -Ax + R(x) + Ay - R(y), x-y\rangle & 
        \leq \|A\| \|x-y\|^2 + \mu \|x-y\|^2 \\
        &= \mu^{*} \|x-y\|^2.
    \end{align*}        
\end{proof}

According to \cite{r7}, the spatially discretized problem \eqref{spatially_discr} is mean-square dissipative if 
$$\alpha^*=2\mu^*+L<0.$$

By the following theorem, we discuss how the mean-square dissipative character of \eqref{spatially_discr} is inherited along the numerical dynamics generated by $\theta$-Maruyama method \eqref{teta}. We point out that $\theta$-IMEX  methods are not taken into consideration in this aspect since, by the analysis of the previous section, we can see that they are not as effective in treating problems of type \eqref{generalSPDE}. 




\medskip

\begin{theorem}
Under Assumptions 3 and 4, the mean-square deviation between two numerical solutions arising from the application of $\theta$-Maruyama method to \eqref{spatially_discr} satisfies the following estimate 
$$\mathbb{E}\left[\|U_{n+1}-Y_{n+1}\|^2\right]\leq \gamma(\theta,\Delta t) \mathbb{E}\left[\|U_{n}-Y_{n}\|^2\right],$$
with 
$$\gamma(\theta,\Delta t)=\frac{1+(1-\theta)^2 \Delta t^2(\| A\|^2+M+2\|A\|\mu) + L \Delta t + (1-\theta)\Delta t \mu^*}{1-2\theta \Delta t \mu^*}.$$
\end{theorem}
\begin{proof}
    In order to study contractivity, we apply $\theta$-Maruyama method \eqref{teta} for time discretization of \eqref{spatially_discr} in correspondence of two distinct solutions, with initial values $U_0$ and $Y_0$, i.e.,
    \begin{equation}
        U^{n+1} = U^{n} + (1-\theta) \Delta t (-AU^n) + \theta \Delta t (-AU^{n+1}) + \Delta t R(U^{n}) + g(U^n) \Delta W_{n},
    \end{equation}
    \begin{equation}
        Y^{n+1} = Y^{n} + (1-\theta) \Delta t (-AY^n) + \theta \Delta t (-AY^{n+1}) + \Delta t R(Y^n) + g(Y^n) \Delta W_{n}.
    \end{equation}
    By applying Lemma from \cite{r7}, and denoting $Z^n=U^n-Y^n$, we obtain
    \begin{equation*}
        (1-2\theta \Delta t \mu^*) \| Z^{n+1}\|^2 \leq \|Z^n + (1-\theta) \Delta t [-A Z^n + R(U^n) - R(Y^n)] + (g(U^n)-g(Y^n))\Delta W_n \|^2,
    \end{equation*}
    whose right-hand side can be bounded by
    \begin{equation*}
        \left[1+(1-\theta)^2 \Delta t^2(\| A\|^2+M+2\|A\|\mu) + L \|\Delta W_n\|^2 + (1-\theta)\Delta t \mu^*\right]\|Z^n\|^2 + M(t), 
    \end{equation*}
    where $M(t)$ is a martingale. Passing to expectation leads to the result.
\end{proof}

In order to study monotonicity and find the timestep restriction needed for it to be preserved, we have to analyze the inequality $\gamma(\theta,\Delta t)<1$, i.e.,
\begin{equation*}
    \frac{1+(1-\theta)^2 \Delta t^2(\| A\|^2+M+2\|A\|\mu) + L \Delta t + (1-\theta)\Delta t \mu^*}{1-2\theta \Delta t \mu^*} < 1.
\end{equation*}
Simplifying, we obtain
\begin{equation}\label{condcontr}
    (1-\theta)^2 \Delta t (\|A\|^2 + M + 2\|A\|\mu) + L + (1-\theta) \mu^* < -2\theta \mu^*,
\end{equation}
which has to be evaluated separately for each case.

For example, when $\theta=1$, we are left with condition
\begin{equation*}
    L + 2\mu^* < 0,
\end{equation*}
which is the contractivity condition of the continuous problem, pointing out that stochastic backward Euler method is unconditionally contractive.

\section{Numerical experiments}
We now provide experimental confirmations of the theoretical results arising from Section 3, using different equations of reaction-diffusion type as examples. As regards Section 4, since the results are mostly highlighting negative behaviours of $\theta$-IMEX methods, their performances are not reported.

Discretisation in space is done by centered finite differences of second order. Higher order discretisations are tested, however they did not yield significantly improved results and, for this reason, they are not reported in the remainder. Accordingly, matrix $A$ is defined by $A=\frac{1}{\Delta x^p}B$, $p=2$ and $B$ is a corresponding Toeplitz matrix, obtained by adding Dirichlet boundary conditions. Upon discretisation in space, discretisation in time is performed using $\theta$-Maruyama method \eqref{teta}. All figures are carried out in semilogarithmic scale.

We observe that contractive condition \eqref{condcontr}, as treated in the previous, is satisfied in the examples where when the noise is additive or the function $g$ is linear. 

\subsection{Ginzburg-Landau equation}
We consider the following Ginzburg-Landau equation \cite{glpaper}:
\begin{equation}\label{GL}
    {\rm d}u=\left(\Delta u+u-u^3\right){\rm d}t +g(u) {\rm d} W,
\end{equation}
covered by class \eqref{generalSPDE}, with operator $\mathcal{A}=-\Delta$ and reaction term ${\cal R}(u)=u-u^3$. 

\begin{figure}
     \centering
     \begin{subfigure}[b]{1\textwidth}
         \centering
         \scalebox{0.5}{\includegraphics{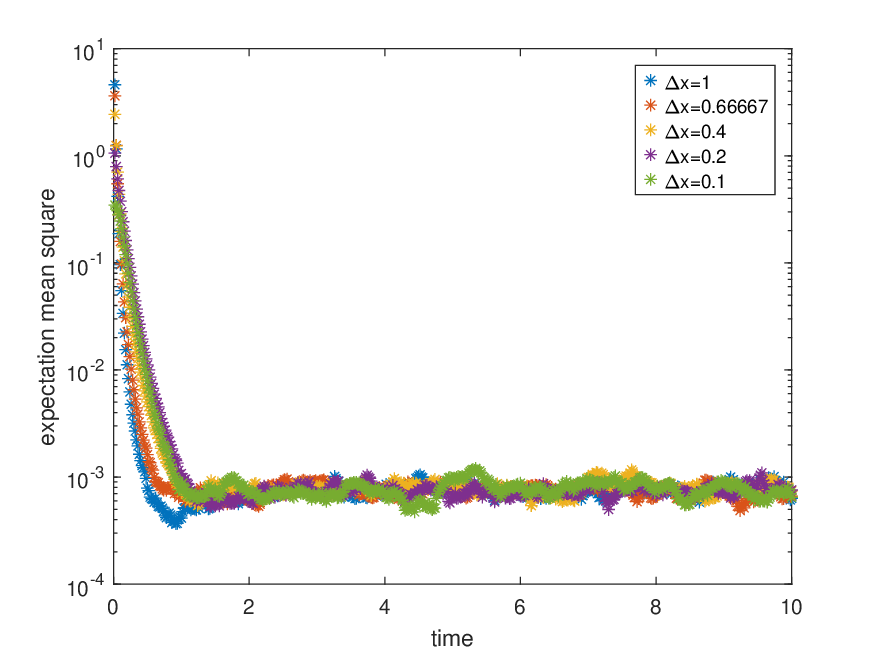}}
         \caption{$\epsilon=0.1$}
         \label{fig:GL_eps01}
     \end{subfigure}
     \vfill
     \begin{subfigure}[b]{1\textwidth}
         \centering
         \scalebox{0.5}
         {\includegraphics{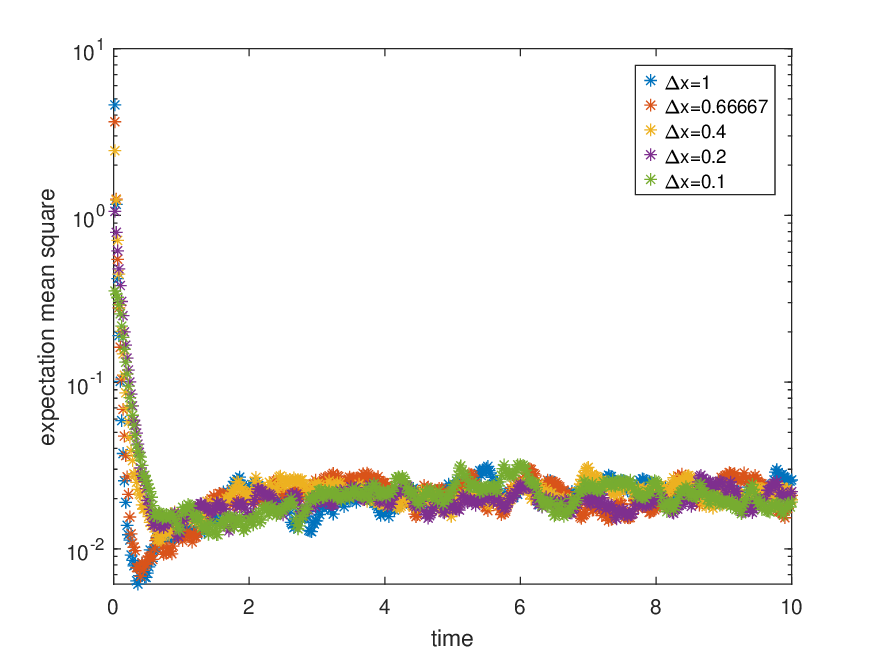}}
         \caption{$\epsilon=0.5$}
         \label{fig:GL_eps05}
     \end{subfigure}
     \vfill
     \begin{subfigure}[b]{1\textwidth}
         \centering
         \scalebox{0.5}
         {\includegraphics{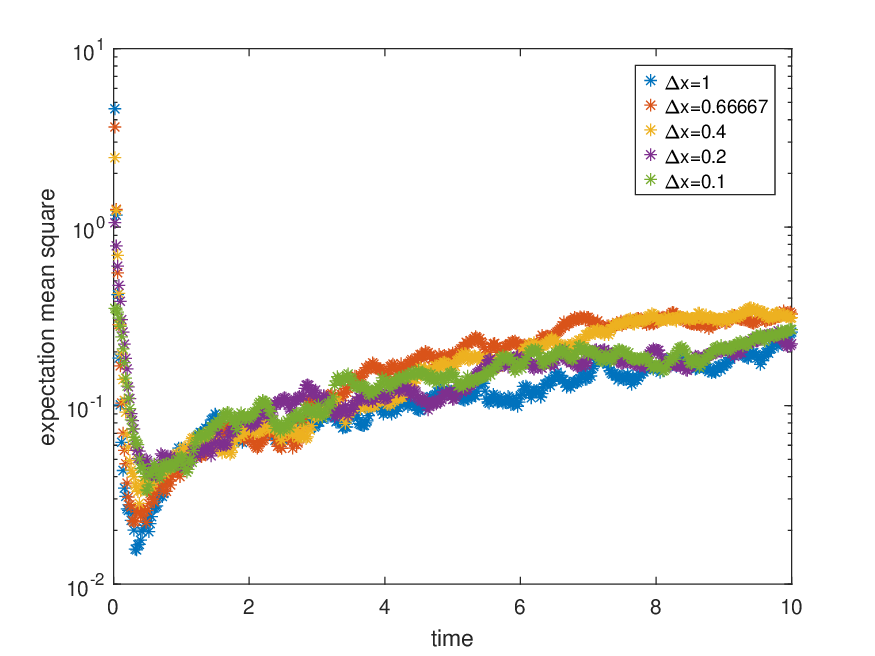}}
         \caption{$\epsilon=0.9$}
         \label{fig:GL_eps09}
     \end{subfigure}
        \caption{Mean-square stability of Ginzburg-Landau equation \eqref{GL} with additive noise for different values of $\epsilon$ using stochastic backward Euler method ($\theta=1$), using 500 steps in time and for different number of points in space.}
        \label{fig:GL_additive}
\end{figure}

In Figure \ref{fig:GL_additive}, the results regarding additive noise are represented, using the stochastic backward Euler method. The noise is regulated by a parameter $g(u)=\epsilon$ and the results are represented for its different values. We can observe that the behaviour of expectation mean-square is similar for different values of $\Delta x$ which confirms the contractivity analysis where we either have contractivity preserved or not without a strong bound on the timestep. Let us observe that, for increasing values of the parameter $\epsilon$, the error in time increases accordingly. This behaviour can be explained in rigorous way, by means of $\epsilon$-expansion arguments \cite{dambook,paternoster}: indeed, one can prove that the numerical methods insert in the corresponding approximate dynamics secular terms growing as $\epsilon\sqrt{t}$. In other terms, the larger the diffusion coefficient is (and, correspondingly, the more dominant the stochastic part of the problem is) or the larger the time window is, the bigger is the contribution of secular terms, affecting the overall accuracy of the numerical solver. 

\begin{figure}
     \centering
     \begin{subfigure}[b]{1\textwidth}
         \centering
         \scalebox{0.5}{\includegraphics{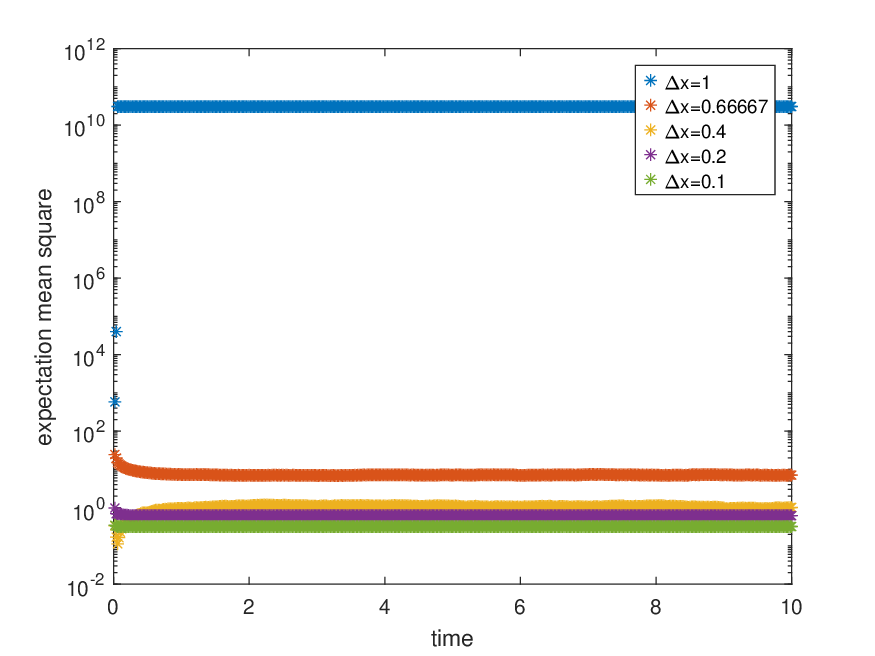}}
         \caption{linear diffusion for $\theta=0$}
         \label{fig:GLa}
     \end{subfigure}
     \vfill
     \begin{subfigure}[b]{1\textwidth}
         \centering
         \scalebox{0.5}
         {\includegraphics{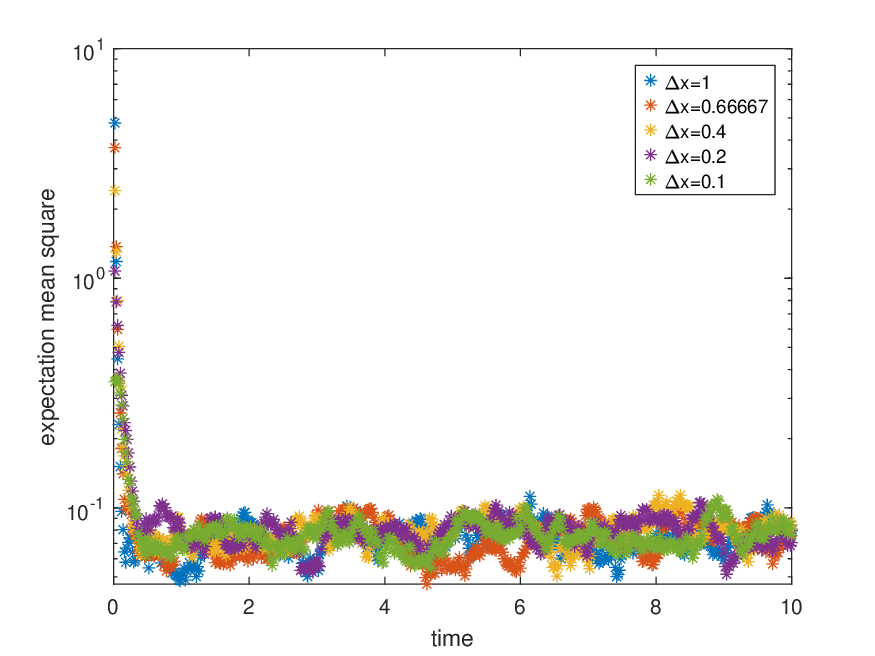}}
         \caption{linear diffusion for $\theta=1$}
         \label{fig:GLb}
     \end{subfigure}
     \vfill
     \begin{subfigure}[b]{1\textwidth}
         \centering
         \scalebox{0.5}
         {\includegraphics{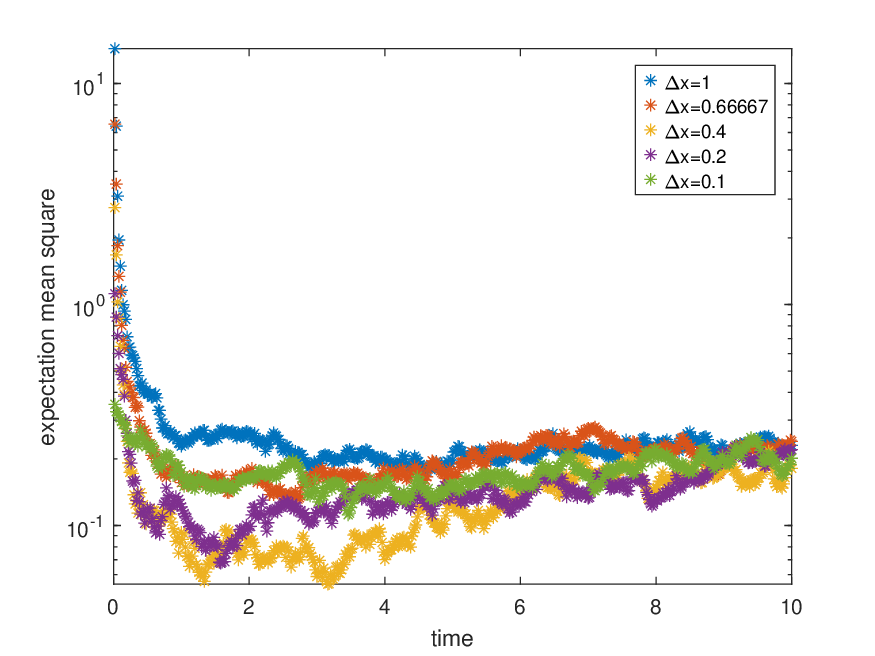}}
         \caption{quadratic diffusion for $\theta=1$}
     \end{subfigure}
        \caption{Mean-square stability of Ginzburg-Landau equation \eqref{GL} with multiplicative noise using Euler-Maruyama method ($\theta=0$) for linear diffusion coefficient and stochastic backward Euler method ($\theta=1$) for both linear and quadratic diffusion coefficient, using 500 points in time and for different number of points in space.}
        \label{fig:GL_multiplicative}
\end{figure}

In Figure \ref{fig:GL_multiplicative}, we present the results for multiplicative noise. In Subfigures \ref{fig:GLa} and \ref{fig:GLb}, the diffusion coefficient is linear, i.e., $g(u)=u$, and tested for both Euler-Maruyama method ($\theta=0$) and stochastic backward Euler method ($\theta=1$). For the Euler-Maruyama method, since $\theta<\frac12$, we can notice that the dissipativity is not preserved. On the other side, for stochastic backward Euler method, we can observe that the behaviour of mean-square expectation is similar for different values of $\Delta x$ which confirms the contractivity analysis. The last subfigure shows the results for quadratic diffusion coefficient, i.e., $g(u)=u^2$. The preservation of dissipativity can be seen for all values of $\Delta x$, as expected from the analysis. 


\subsection{Cahn-Hilliard equation}
We consider the following Cahn-Hilliard equation \cite{ch}:
\begin{equation}\label{CH}
    {\rm d}u=\left(-\Delta^2 u + \Delta V'(u)\right){\rm d}t +g(u) {\rm d} W,
\end{equation}
where $V(u(x,t))$ is the potential and, in correspondence with \eqref{generalSPDE}, the operator $\mathcal{A}=\Delta^2$ and the reaction term ${\cal R}(u)=\Delta V'(u)$. As a potential, we choose $V(u)=u-u^3$, as often found in literature. However, the results presented below were similar also for different choices of the potential.

\begin{figure}
        \begin{subfigure}{.5\textwidth}
            \scalebox{0.45}{\includegraphics{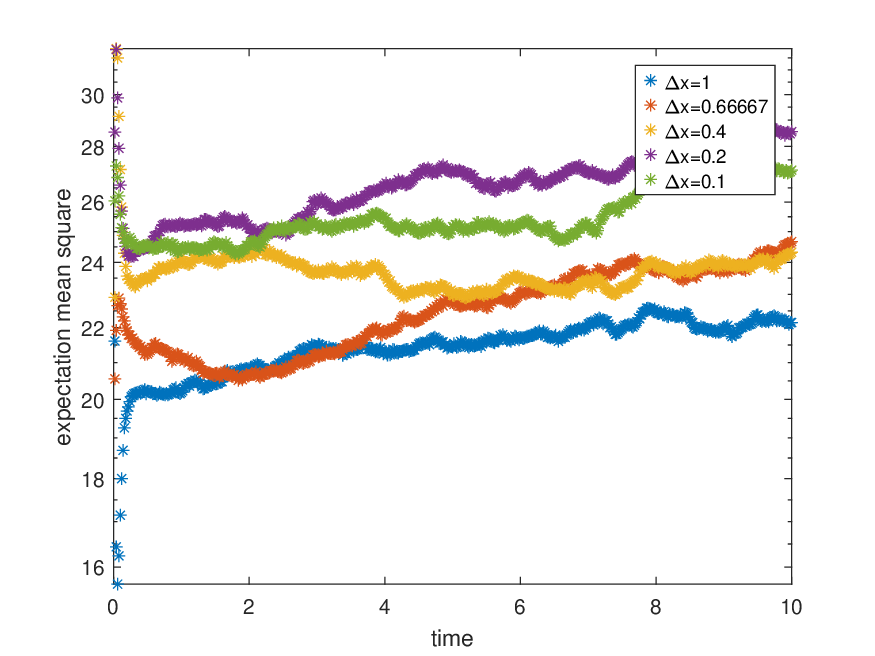}}
            \subcaption{additive noise}
        \end{subfigure}
        \hfill
        \begin{subfigure}{.5\textwidth}
            \scalebox{0.45}{\includegraphics{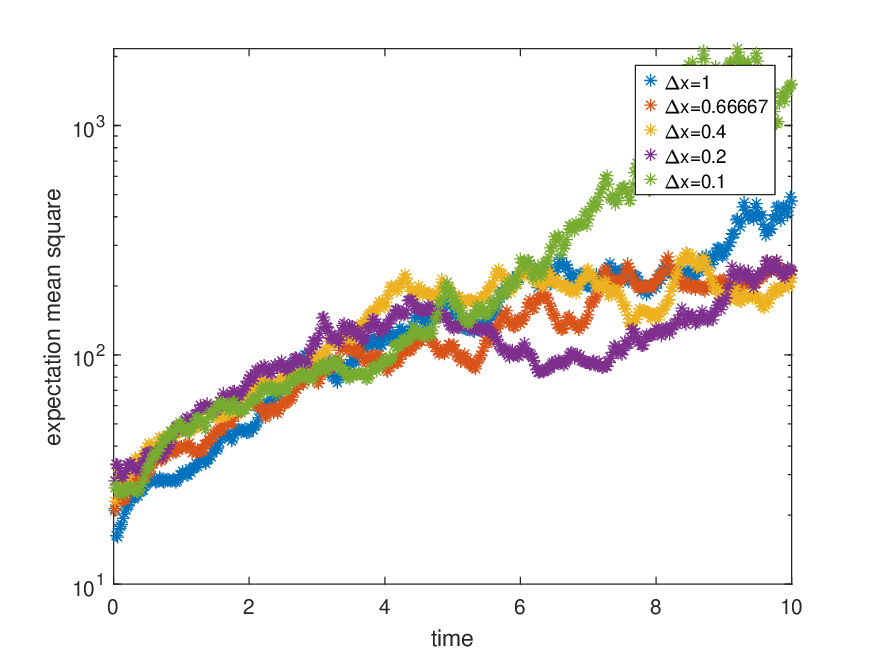}}
            \subcaption{linear multiplicative noise}
        \end{subfigure}
    \caption{Numerical results for Cahn-Hilliard problem \eqref{CH}}
    \label{fig:Cahn-Hilliard}
\end{figure}

In Figure \ref{fig:Cahn-Hilliard}, we present results on dissipativity behaviour of stochastic backward Euler method, for both additive and linear multiplicative noises. Even if we are aware of the dissipative character of this problem \cite{weinan}, experiments do not confirm its preservation. Additionally, in the case of linear multiplicative noise, we can observe that the exponential mean-square of the difference between two solutions is blowing up. The reasoning for this behaviour in both cases comes from another important observation in the contractivity analysis, that is the dependence on the relation between $\Delta t$ and $\Delta x$ that drastically impacts the contractivity conservation. 
By Toeplitz discretization, we have that
\begin{equation*}
    \| \Delta^2 u \| \sim \frac{1}{\Delta x^4}
\end{equation*}
which imposes a stronger constraint on relation between $\Delta t$ and $\Delta x$ such that
\begin{equation*}
    \Delta t < \frac{\Delta x^8}{2}
\end{equation*}

This example is illustrating exactly when we should be cautious while dealing with stochastic reaction-diffusion equations in the context of dissipativity preservation by stochastic $\theta$-methods. In conclusion, even if this example is a negative one in terms of dissipativity preservation, it provides a relevant test useful to highlight the importance of keeping the ratio between time and spatial stepsize bounded. In some cases, as it is for Cahn-Hilliard problem, this restriction appears to be highly severe. 

\subsection{Example 3. Uncoupled system}

We consider the following uncoupled system of SPDEs 
\begin{equation}\label{uncoupled}
    {\rm d} \begin{pmatrix} u \\ v \end{pmatrix} = \begin{pmatrix} \Delta & 0 \\ 0 & \Delta \end{pmatrix} \begin{pmatrix} u \\ v \end{pmatrix} + \begin{pmatrix}u-u^3 \\ v \end{pmatrix} + \begin{pmatrix} g_1(u) \\ g_2(u) \end{pmatrix} {\rm d} W,
\end{equation}
arising from coupling a Ginzburg-Landau reaction term with a linear one. It is covered by the general class \eqref{generalSPDE}, with operator $\mathcal{A}=-\Delta$, reaction term $${\cal R}\left(\begin{pmatrix} u \\ v \end{pmatrix}\right)=\begin{pmatrix}u-u^3 \\ v \end{pmatrix},$$
and $g_1(u),g_2(u)$ being components of the diffusion coefficient $g(u)$.
 
\begin{figure}
     \centering
     \begin{subfigure}[b]{1\textwidth}
         \centering
         \scalebox{0.5}{\includegraphics{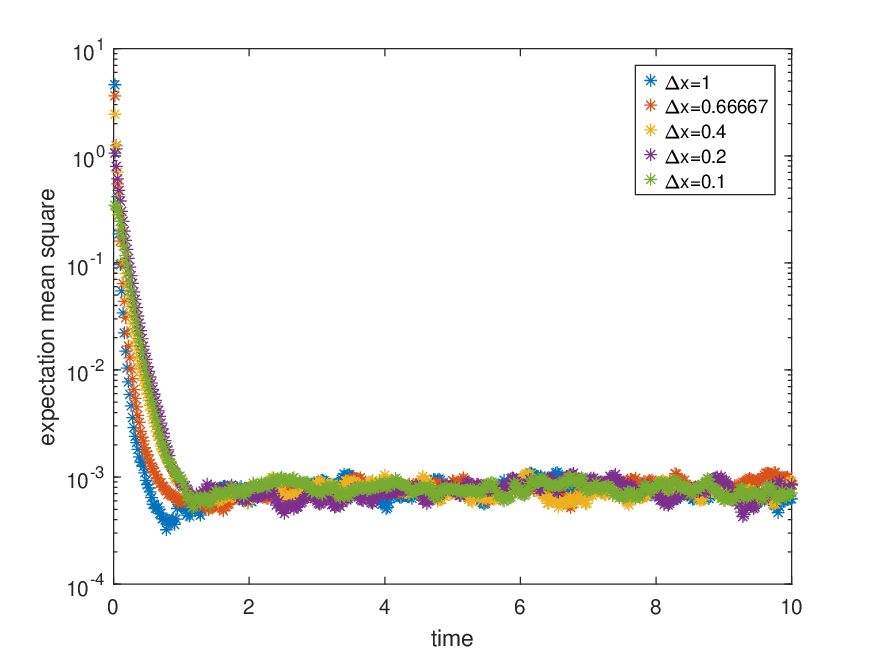}}
         \caption{constant diffusion for $\theta=1$}
     \end{subfigure}
     \vfill
     \begin{subfigure}[b]{1\textwidth}
         \centering
         \scalebox{0.5}
         {\includegraphics{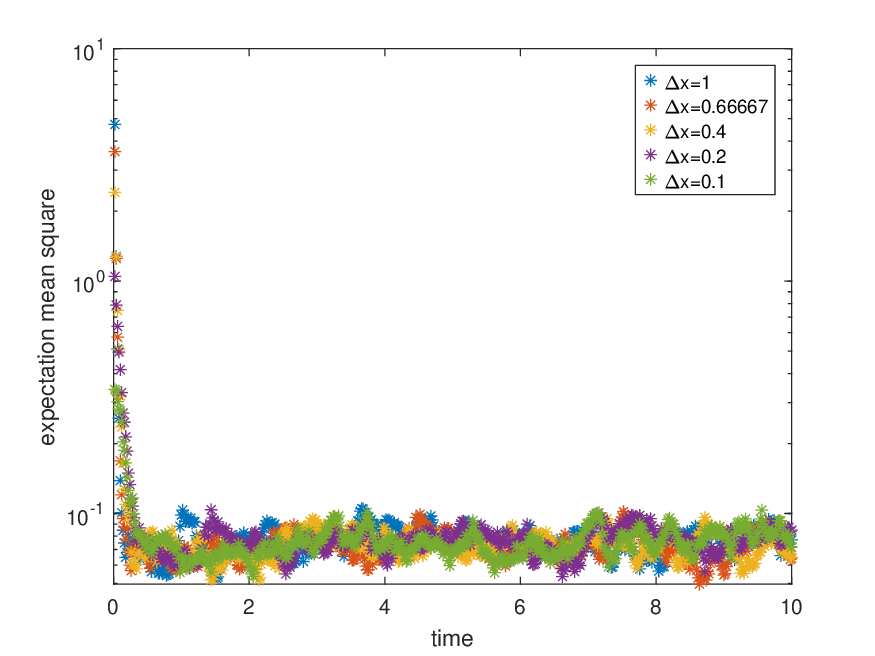}}
         \caption{linear diffusion for $\theta=1$}
     \end{subfigure}
     \vfill
     \begin{subfigure}[b]{1\textwidth}
         \centering
         \scalebox{0.5}
         {\includegraphics{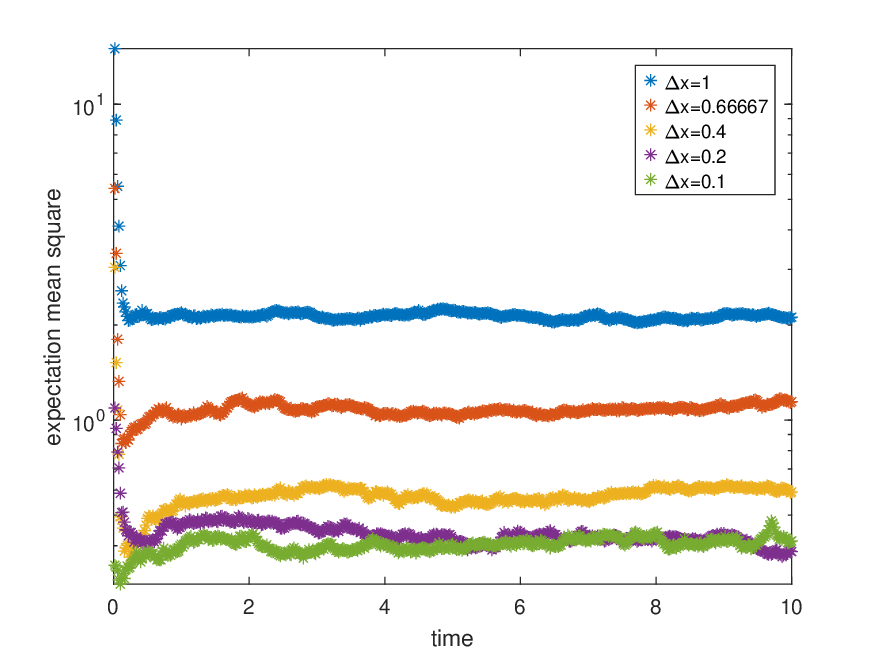}}
         \caption{quadratic diffusion for $\theta=1$}
         \label{fig:uncoupled_QU}
     \end{subfigure}
        \caption{Numerical results for the uncoupled system of SPDEs \eqref{uncoupled}}
        \label{fig:uncoupled}
\end{figure}

In Figure \ref{fig:uncoupled}, we present the results for the uncoupled system of SPDEs, using stochastic backward Euler method ($\theta=1$), for both additive noise and multiplicative linear and quadratic diffusion coefficients. For additive and multiplicative linear noises, a similar behaviour of exponential mean-square contractivity can be observed for different values of $\Delta x$, which is in concordance with the contractivity analysis. This result is reinforcing our theory to be applied even for systems of SPDEs in feasible way.

It is worth highlighting that, in case of systems of SPDEs with quadratic diffusion coefficient, a spatial stepsize restriction for the conservation of dissipativity looks evident in Subfigure \ref{fig:uncoupled_QU}. This behaviour reasonably matches with the presence of $L_g^2\Delta t$ in the numerator of $\alpha(\theta,\Delta t)$ in \eqref{coeff2}: indeed, when $\Delta t$ is fixed and the Lipschitz constant is large enough, a significant counterpart of the spatial stepsize to numerically keep the contractive behaviour is expected.

\subsection{Example 4. DIB system}

We consider the following stochastic perturbation of DIB system \cite{DIB}:
\begin{equation}\label{DIB}
\begin{aligned}
    {\rm d} \begin{pmatrix} u \\ v \end{pmatrix} &= \begin{pmatrix} d_1 \Delta & 0 \\ 0 & d_2 \Delta \end{pmatrix} \begin{pmatrix} u \\ v \end{pmatrix} \\[2mm]
    &+ \rho \begin{pmatrix} A_1 (1-v) u - A_2 u^3 - B(v - \alpha) \\ C(1+k_2 u)(1-v)(1-\gamma(1-v)) - D v(1+k_3u) (1+\gamma v) \end{pmatrix} + \begin{pmatrix} g_1(u) \\ g_2(u) \end{pmatrix} {\rm d} W
    \end{aligned}
\end{equation}
covered by class \eqref{generalSPDE}, with operator $\mathcal{A}=-\Delta$, reaction term $${\cal R}\left(\begin{pmatrix} u \\ v \end{pmatrix}\right)=\rho \begin{pmatrix} A_1 (1-v) u - A_2 u^3 - B(v - \alpha) \\ C(1+k_2 u)(1-v)(1-\gamma(1-v)) - D v(1+k_3u) (1+\gamma v) \end{pmatrix},$$
and $g_1(u),g_2(u)$ being components of the diffusion coefficient $g(u)$.

This system is describing an electrodeposition process, for instance in battery modelling, and represents one of the cases of Turing patterns formation.

\begin{figure}
     \centering
     \begin{subfigure}[b]{1\textwidth}
         \centering
         \scalebox{0.5}{\includegraphics{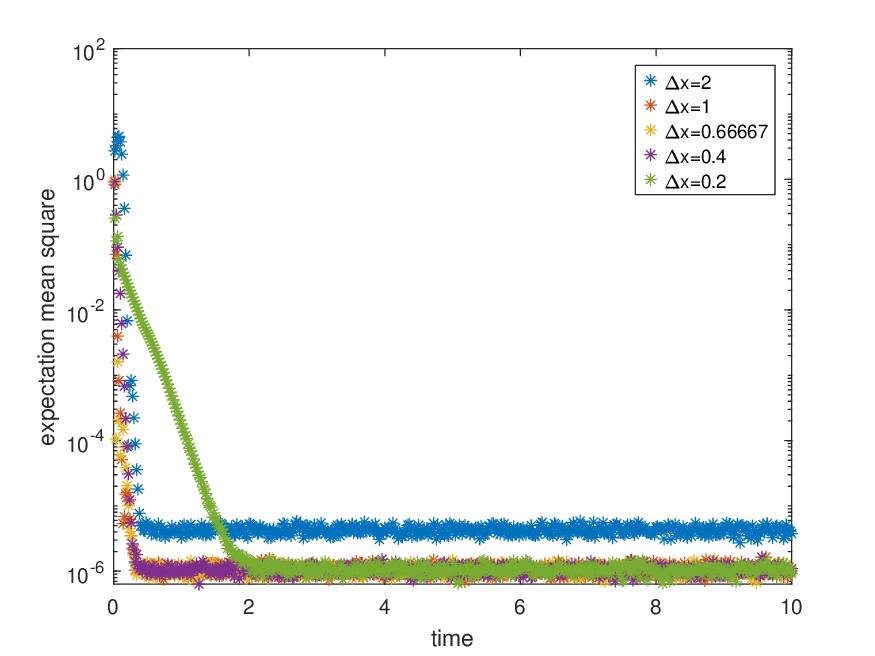}}
         \caption{constant diffusion for $\theta=1$}
     \end{subfigure}
     \vfill
     \begin{subfigure}[b]{1\textwidth}
         \centering
         \scalebox{0.5}
         {\includegraphics{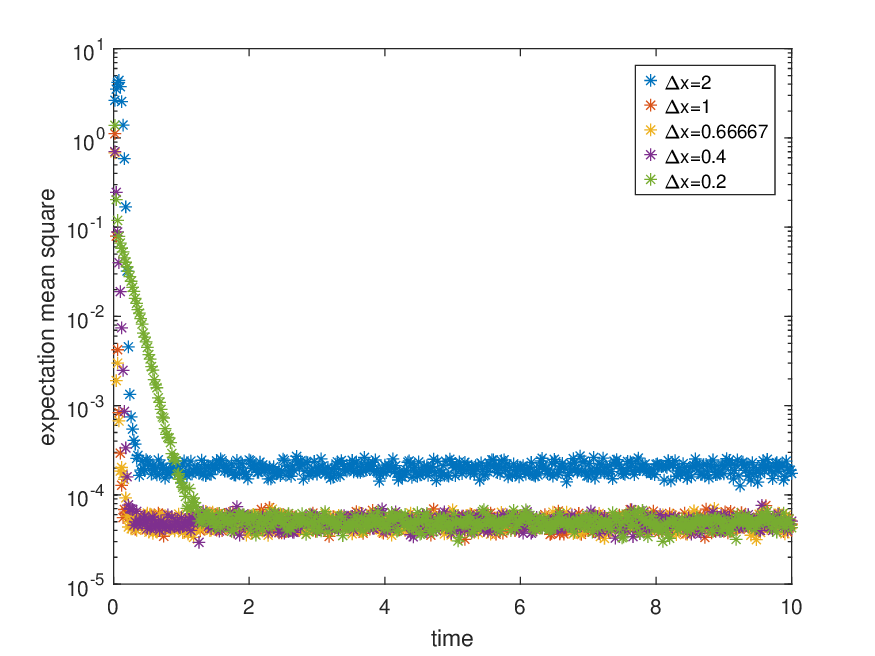}}
         \caption{linear diffusion for $\theta=1$}
     \end{subfigure}
     \vfill
     \begin{subfigure}[b]{1\textwidth}
         \centering
         \scalebox{0.5}
         {\includegraphics{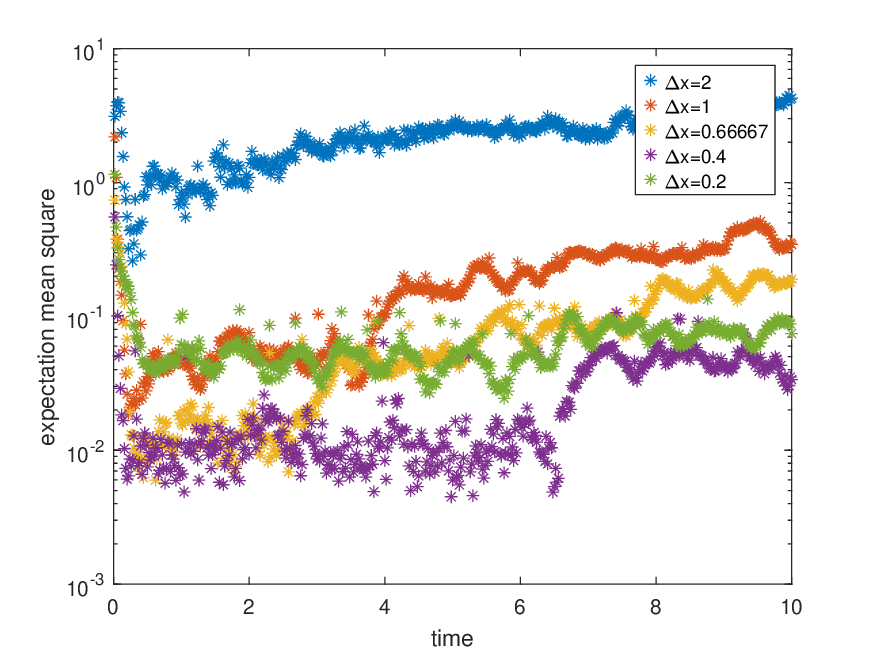}}
         \caption{quadratic diffusion for $\theta=1$}
     \end{subfigure}
        \caption{Numerical results for DIB system \eqref{DIB}}
        \label{fig:DIB}
\end{figure}

In Figure \ref{fig:DIB}, we present the results for the DIB system of SPDEs, using stochastic backward Euler method ($\theta=1$), for both additive noise and multiplicative linear and quadratic diffusion coefficient. We can observe a similar behaviour of mean-square dissipativity conservation as before, for different value of $\Delta x$, in concordance with the contractivity analysis. We draw attention to case of quadratic diffusion, where the choice of the spatial and temporal stepsize plays an important role as in previous example. 

\section{Conclusions and open issues} 
The paper has handled the numerical solution of stochastic reaction-diffusion problems, giving rise to a theoretical framework for the analysis of the conservation of mean-square dissipativity in the time integration of the spatially discretized problem via finite differences. We have considered both stochastic $\theta$-methods and stochastic IMEX methods, under the assumptions of local Lipschitz continuity of the diffusion coefficient and $C^{1}$ continuity of the drift, eventually also equipped by its one-sided Lipschitz continuity. The investigation has clarified the conservative character of these methods, in terms of preservation of mean-square contractivity, and it has also been supported by selected numerical experiments.

This work opens the horizons towards answering several open problems regarding stochastic PDEs. The paper that mostly inspired our investigation \cite{weinan} gives a complete theory of invariant measure, as a way to treat dissipativity, exactly for the type of the equations for which we provided numerical theory of mean-square stability. Theoretical results presented there together with our numerical treatment give a sufficient foundation to start working on a numerical treatment of invariant measure, that is a future perspective of this research. Additionally, it is worth observing that a full stability analysis of stochastic IMEX methods needs to be well established in the scientific literature, as well as further aspects of numerical conservation of invariance laws that characterize stochastic PDEs.

\section*{Declarations}
\begin{itemize}
\item {\bf Funding}. The authors are members of the INdAM Research group GNCS. This work is supported by GNCS-INDAM project and by  PRIN-PNRR project P20228C2PP BAT-MEN (BATtery Modeling, Experiments \& Numerics) - Enhancing battery lifetime: mathematical modeling, numerical simulations and AI parameter estimation techniques (CUP: E53D23017940001).
\item {\bf Conflict of interest}. The authors declare no competing interests.
\item {\bf Author contribution}. The authors contributed equally to this work.
\end{itemize}

\end{document}